\documentclass{article}
\usepackage{latexsym,amsfonts,amsmath,amsthm,amssymb,makeidx}
\usepackage[title]{appendix}

\let\oldbibliography\thebibliography
\renewcommand{\thebibliography}[1]{%
\oldbibliography{#1}%
\setlength{\itemsep}{0pt}%
}

\makeindex
\newtheorem{definition}{Definition}[section]
\newtheorem{theorem}{Theorem}[section]
\newtheorem{lemma}{Lemma}[section]
\newtheorem{corollary}{Corollary}[section]

\newtheorem{remark}{Remark}[section]
\newcommand{\RN}{\mathbb R^N}
\newcommand{\Om}{\Omega}

\newcommand{\iy}{\infty}

\newcommand{\dd}{\delta}

\newcommand{\la}{\lambda}

\newcommand{\R}{\mathbb R}

\newcommand{\bb}{\beta}

\newcommand{\e}{\varepsilon}
\newcommand{\vp}{\varphi}

\newcommand{\bt}{\begin{theorem}}
\newcommand{\et}{\end{theorem}}
\newcommand{\bl}{\begin{lemma}}
\newcommand{\el}{\end{lemma}}
\newcommand{\bd}{\begin{definition}}
\newcommand{\ed}{\end{definition}}
\newcommand{\bc}{\begin{corollary}}
\newcommand{\ec}{\end{corollary}}
\newcommand{\bp}{\begin{proof}}
\newcommand{\ep}{\end{proof}}
\newcommand{\bx}{\begin{example}}
\newcommand{\ex}{\end{example}}
\newcommand{\bi}{\begin{exercise}}
\newcommand{\ei}{\end{exercise}}
\newcommand{\bo}{\begin{prop}}
\newcommand{\eo}{\end{prop}}
\newcommand{\br}{\begin{remark}}
\newcommand{\er}{\end{remark}}
\newcommand{\be}{\begin{equation}}
\newcommand{\ee}{\end{equation}}
\newcommand{\ba}{\begin{align}}
\newcommand{\ea}{\end{align}}
\newcommand{\bn}{\begin{enumerate}}
\newcommand{\en}{\end{enumerate}}
\newcommand{\bg}{\begin{align*}}
\newcommand{\bcs}{\begin{cases}}
\newcommand{\ecs}{\end{cases}}

\newcommand{\Sg}{\Sigma}

\newcommand{\sg}{\sigma}
\newcommand{\bean}{\begin{eqnarray*}}
\newcommand{\eean}{\end{eqnarray*}}

\numberwithin{equation}{section}

\begin{document}

\title{\bf Removable singularity of positive solutions for a critical
elliptic system with isolated singularity \thanks{
E-mail address:
chenzhijie1987@sina.com, zjchen@tims.ntu.edu.tw(Chen); cslin@math.ntu.edu.tw(Lin)}}
\date{}
\author{{\bf Zhijie Chen$^1$, Chang-Shou Lin$^2$}\\
\footnotesize {\it  $^{1}$Center for Advanced Study in Theoretical Sciences, National Taiwan
University, Taipei 106, Taiwan}\\
\footnotesize {\it  $^{2}$Taida Institute for Mathematical Sciences,
Center for Advanced Study in}\\
\footnotesize {\it Theoretical Sciences, National Taiwan
University, Taipei 106, Taiwan}
}

\maketitle
\begin{center}
\begin{minipage}{120mm}
\begin{center}{\bf Abstract}\end{center}

We study qualitative properties of positive singular solutions to a two-coupled elliptic system with critical exponents. This system is related to coupled nonlinear Schr\"{o}dinger equations
with critical exponents for nonlinear optics and Bose-Einstein condensates.
We prove a sharp result on the removability of the same isolated singularity for both two components of the solutions.
We also prove the nonexistence of positive solutions with one component bounded near the singularity and the other component unbounded near the singularity. These results will be applied in a subsequent work where the same system in a punctured ball will be studied.

\vskip0.10in


\end{minipage}
\end{center}

\vskip0.10in
\section{Introduction}
\renewcommand{\theequation}{1.\arabic{equation}}

In this paper, we study qualitative properties of positive singular solutions to the following two-coupled elliptic system with critical exponents:
\be\label{eq0}
\begin{cases}
\begin{split}-\Delta u =
\mu_1 u^{2^\ast-1}+\beta u^{\frac{2^\ast}{2}-1} v^{\frac{2^\ast}{2}}\\
-\Delta v =\mu_2 v^{2^\ast-1}+\beta v^{\frac{2^\ast}{2}-1}u^{\frac{2^\ast}{2}}\\
\end{split}  \quad   \text{in $\RN\setminus\{0\}$},\\
u, v> 0 \quad\hbox{and}\quad
u, v\in C^2(\RN\setminus\{0\}).\end{cases}\ee
Here, $\mu_1, \mu_2$ and $\bb$ are all positive constants, $N\ge 3$ and $2^\ast:=\frac{2N}{N-2}$ is the Sobolev critical exponent.

System (\ref{eq0}) is the critical case of the following general problem
\be\label{eq002}
\begin{cases}
\begin{split}-\Delta u =\la_1 u +
\mu_1 u^{2p-1}+\beta u^{p-1}v^{p}\\
-\Delta v =\la_2 v+\mu_2 v^{2p-1}+\beta v^{p-1} u^{p}\\
\end{split}  \quad   \text{in $\Om$},\\
u, v> 0 \,\,\hbox{in $\Om$}\quad\hbox{and}\quad
u, v\in C^2(\Om),\end{cases}\ee
where $\Om\subset\RN$ is a domain, $p>1$ and $2p\le 2^\ast$ if $N\ge 3$.
In the case $p=2$, the cubic system (\ref{eq002}) arises in mathematical models from various physical phenomena, such as the incoherent solitons in
nonlinear optics and the multispecies Bose-Einstein condensates in hyperfine spin states. We refer the reader for these to the survey articles \cite{F, KL} which, among other things,
also contain information about the physical relevance of non-cubic nonlinearities.
In the subcritical case $2p<2^\ast$, this problem has
been widely investigated via variational methods and
topological methods in the last decades; see \cite{CLLL, CZ3, LW1, MMP, QS, S}
and references therein.

On the other hand, critical exponent problems, including various scalar equations and coupled systems, are very interesting and challenging in view of mathematics, and have received ever-increasing interest and have been studied intensively in the literature. For system (\ref{eq002}) in the critical case $2p=2^\ast$, we refer the reader to \cite{ChenLin, CLZ, CZ, CZ2, SK} for the bounded domain case, where the existence and asymptotic behaviors of positive least energy solutions and sign-changing solutions were well investigated. For the entire space case $\Om=\RN$, if $\la_1, \la_2>0$, it is easy to see that (\ref{eq002}) has no solutions via the method of moving planes. Therefore, in the entire space case $\Om=\RN$ with $2p=2^\ast$, we always assume $\la_1=\la_2=0$. Then system (\ref{eq002}) turns to be
\be\label{eq1}
\begin{cases}
\begin{split}-\Delta u =
\mu_1 u^{2^\ast-1}+\beta u^{\frac{2^\ast}{2}-1} v^{\frac{2^\ast}{2}}\\
-\Delta v =\mu_2 v^{2^\ast-1}+\beta v^{\frac{2^\ast}{2}-1}u^{\frac{2^\ast}{2}}\\
\end{split}  \quad   \text{in $\RN$},\\
u, v> 0 \quad\hbox{and}\quad
u, v\in C^2(\RN),\end{cases}\ee
which has also been well studied in \cite{ChenLi, GLW, GL, LM}. See also \cite{AFP, CZ-CVPDE2013, CZ-TAMS2013} for the study of a more general system (i.e. system (\ref{eq1}) is a special case of it) with critical exponents and Hardy potentials. Remark that system (\ref{eq1}) relies essentially on the sign of the coupling constant $\bb$. For the replusive case $\bb<0$ that we do not consider in this paper, Guo, Li and Wei \cite{GLW} proved the existence of infinitely many {\it non-radial} solutions to (\ref{eq1}) for the low dimensional case $N\le 4$. While for the case $\bb>0$ as considered here, it is natural to suspect that every solution of (\ref{eq1}) is {\it radially symmetric} up to a translation, which has been proved by Guo and Liu \cite{GL} and Chen and Li \cite{ChenLi} independently. In fact, they studied more general systems and (\ref{eq1}) is a special case of their problems.\\

\noindent {\bf Theorem  A.} (\cite{ChenLi, GL})
{\it Assume that $\mu_1, \mu_2, \bb>0$ and $(u, v)$ is an entire solution of (\ref{eq1}). Then $u$ and $v$ are radially symmetric with respect to the same point. In particular, $(u, v)=(k U, l U)$, where $k, l>0$ satisfies
\be\label{eq3}\mu_1 k^{2^\ast-2}+\bb k^{\frac{2^\ast}{2}-2}l^{\frac{2^\ast}{2}}=1,\quad \mu_2 l^{2^\ast-2}+\bb l^{\frac{2^\ast}{2}-2}k^{\frac{2^\ast}{2}}=1,\ee
and $U>0$ is an entire solution of
\be\label{eq2}-\Delta u=u^{2^\ast-1},\quad u>0\,\,\,\text{in $\RN$},\ee
namely, there exists $\epsilon>0$ and $y\in\RN$ such that
$$U(x)=[N(N -2)]^{\frac{N-2}{4}}\left(\frac{\epsilon}{\epsilon^2 +
|x-y|^2}\right)^{\frac{N-2}{2}}.$$ }

Theorem A gives a complete classification of entire solutions to system (\ref{eq1}). As far as we know, there seems {\it no} any result about positive singular solutions of problem (\ref{eq1}). This is one of our motivations to study problem (\ref{eq0}), a solution of which has at most one isolated singularity and is the simplest example of positive singular solutions to problem (\ref{eq1}). Another motivation for studying problem (\ref{eq0}) comes from the study of positive singular solutions of the conformal scalar curvature equation (\ref{eq2}).
In view of conformal geometry, a solution $u$ of (\ref{eq2}) defines a conformally flat metric $g_{ij}=u^{\frac{4}{N-2}}\dd_{ij}$ with constant scalar curvature.
The classical work by Schoen and Yau \cite{Schoen, Schoen1, SchoenY} on conformally flat manifolds and the Yamabe problem has highlighted the importance of studying solutions of (\ref{eq2}) with a nonempty singular set. The simplest case of (\ref{eq2}) with a nonempty singular set is the following problem
\be\label{eq4}
\begin{cases}
-\Delta u = u^{2^\ast-1}  \quad   \text{in $\RN\setminus\{0\}$},\\
u> 0 \quad\hbox{and}\quad
u\in C^2(\RN\setminus\{0\}).\end{cases}\ee
Problem (\ref{eq4}) and related problems have received great interest and have been widely studied in \cite{CGS, CL1, CL2, CL3, Fowler, KMPS, Lin, TZ} and references therein. In particular, Fowler \cite{Fowler} described all radial solutions of (\ref{eq4}), and Caffarelli, Gidas and Spruck \cite{CGS} proved the radial symmetry of all solutions of (\ref{eq4}). Thus, all solutions of (\ref{eq4}) can be well understood through the study of the corresponding ordinary differential equation. Denote $B_r:=\{x\in\RN : |x|<r\}$. By the Pohozaev identity, we have for $0<s<r$ that
$P(r; u)-P(s; u)=0$,
where
\be\label{eq5}P(r; u)=\int_{\partial B_r}\left(\frac{N-2}{2}u\frac{\partial u}{\partial \nu}-\frac{r}{2}|\nabla u|^2+r\left|\frac{\partial u}{\partial \nu}\right|^2+\frac{r}{2^\ast}u^{2^\ast}\right)d\sigma,\ee
and $\nu$ is the unit outer normal of $\partial B_r$. Hence $P(r; u)$ is a constant independent of $r$, and we denote it by $P(u)$. Then a classical result of removable singularity is as follows.\\

\noindent {\bf Theorem  B.} (\cite{CGS, Fowler})
{\it Let $u$ be a solution of (\ref{eq4}), then $P(u)\le 0$. Moreover, $P(u)=0$ if and only if $0$ is a removable singularity of $u$, namely $u$ is an entire solution of (\ref{eq2}).}\\

Later Chen and the second author \cite{CL3} extended Theorem B to a more general problem $-\Delta u=K(x)u^{2^\ast-1}$, where $K(x)$ is a positive function satisfying some further conditions; see \cite[Theorem 1.1]{CL3}. Furthermore, they also give conditions on $K(x)$ such that the counterpart of Theorem B does not hold; see \cite[Theorem 1.6]{CL3} for details.

Observe that when $\bb=0$, problem (\ref{eq0}) is just the problem (\ref{eq4}) up to a multiplication.
Similarly as (\ref{eq4}), we can also treat problem (\ref{eq0}) in view of conformal geometry.
A solution $(u, v)$ of (\ref{eq0}) defines two conformally flat metric $g^1_{ij}=u^{\frac{4}{N-2}}\dd_{ij}$ and $g^2_{ij}=v^{\frac{4}{N-2}}\dd_{ij}$. The metrics $g^1$ and $g^2$ pertain to positive scalar curvature functions \be\label{curvature}\mathcal{K}_1=\frac{4(N-1)}{N-2}\left(\mu_1+\bb \left(\frac{v}{u}\right)^{\frac{2^\ast}{2}}\right) \quad\text{and}\quad \mathcal{K}_2=\frac{4(N-1)}{N-2}\left(\mu_2+\bb \left(\frac{u}{v}\right)^{\frac{2^\ast}{2}}\right)\ee
respectively, namely the scalar curvature functions are {\it not prescribed} but are determined by the solution $(u, v)$ itself. Remark that $(\frac{N-2}{4(N-1)}\mathcal{K}_1-\mu_1)(\frac{N-2}{4(N-1)}\mathcal{K}_2-\mu_2)\equiv\bb^2$. An interesting question is {\it whether $\mathcal{K}_1$ and $\mathcal{K}_2$ defined by $(u, v)$ in (\ref{curvature}) are necessarily constants or not}. Clearly it is equivalent to a basic question {\it whether a solution $(u, v)$ of (\ref{eq0}) satisfies $u/v\equiv\text{constant} $ or not}, which is also interesting in itself in view of Theorem A. See also a recent work \cite{MSS} for other systems on the property of proportionality of components.

The main purpose of this paper is to extend Theorem B to system (\ref{eq0}). As a corollary, we can give a partial answer to the above question. For any solution $(u, v)$ of system (\ref{eq0}), we can prove that $\lim_{|x|\to 0}u(x)$ and $\lim_{|x|\to 0} v(x)$ make sense (the limit may be $+\iy$ in general); see Theorem \ref{th2-1} in Section 2. Hence there are three possibilities in general: (1) both $\lim_{|x|\to 0}u(x)<+\iy$ and $\lim_{|x|\to 0}v(x)<+\iy$, namely $(u, v)$ is an entire solution of (\ref{eq1}); (2) $\lim_{|x|\to 0}u(x)=\lim_{|x|\to 0}v(x)=+\iy$, namely the origin is a non-removable singularity of both $u$ and $v$, and we call solutions of this type {\it both-singular solutions}; (3) either $\lim_{|x|\to 0}u(x)=+\iy>\lim_{|x|\to 0}v(x)$ or $\lim_{|x|\to 0}v(x)=+\iy>\lim_{|x|\to 0}u(x)$, and we call solutions of this type {\it semi-singular solutions}.

Multiplying the first equation of (\ref{eq0}) by $x\cdot\nabla u$, the second equation by $x\cdot\nabla v$, and integrating over $B_r\setminus B_s$, we easily obtain the following Pohozaev identity
\be\label{eq6}K(r; u, v)=K(s; u, v),\ee
where
{\allowdisplaybreaks
\begin{align}\label{eq7}
&K(r; u, v)=\int_{\partial B_r}\Bigg[\frac{N-2}{2}\left(u\frac{\partial u}{\partial \nu}+v\frac{\partial v}{\partial \nu}\right)-\frac{r}{2}(|\nabla u|^2+|\nabla v|^2)\nonumber\\
&\quad+r\left|\frac{\partial u}{\partial \nu}\right|^2+r\left|\frac{\partial v}{\partial \nu}\right|^2+\frac{r}{2^\ast}\left(\mu_1 u^{2^\ast}+\mu_2 v^{2^\ast}+2\bb u^{\frac{2^\ast}{2}}v^{\frac{2^\ast}{2}}\right)\Bigg]d\sigma.
\end{align}
}%
Hence $K(r; u, v)$ is a constant independent of $r$, and we denote this constant by $K(u, v)$.
Our main result is about removable singularity, which extends the classical result Theorem B to system (\ref{eq0}).

\bt\label{th1}Let $\mu_1, \mu_2, \bb>0$ and $(u, v)$ be a solution of (\ref{eq0}). Then $K(u, v)\le 0$.
Furthermore, $K(u, v)=0$ if and only if $u, v\in C^2(\RN)$, namely both $u$ and $v$ are smooth at $0$.\et

As an application of Theorems A and \ref{th1}, we obtain the following result immediately.

\begin{corollary}\label{corollary1-1}
Let $\mu_1, \mu_2, \bb>0$ and $(u, v)$ be a solution of (\ref{eq0}) with $K(u, v)=0$, then $(u, v)=(kU, lU)$, where $(k, l)$ and $U$ are seen in Theorem A. In particular, the scalar curvature functions $\mathcal{K}_1$ and $\mathcal{K}_2$ defined by $(u, v)$ in (\ref{curvature}) are constants.
\end{corollary}

Theorem 1.1 indicates that system (\ref{eq0}) has no solutions satisfying $K(u, v)>0$. Moreover, the origin is a removable singularity of both $u$ and $v$ and so $(u, v)$ satisfies Theorem A, provided $K(u, v)=0$. In other words, $(u, v)$ is a singular solution (either both-singular or semi-singular) if and only if $K(u, v)<0$. To the best of our knowledge, this seems to be the {\it first} result on the removability of the same isolated singularity for this kind of critical systems.

Our second result is concerned with singular solutions. A basic question is {\it whether semi-singular solutions exist or not}. Remark that, if $(u, v)$ is a semi-singular solution of (\ref{eq0}), then $u/v\not\equiv\text{constant}$ and so neither $\mathcal{K}_1$ nor $\mathcal{K}_2$ is a constant. Here we can only give a partial answer.

\bt\label{th2}Let $\mu_1, \mu_2, \bb>0$ and $N\ge 3$.
\begin{itemize}
\item[$(1)$]
System (\ref{eq0}) has no semi-singular solutions whenever $N\ge 4$.

\item[$(2)$]
Assume in addition that, if $N=3$ system (\ref{eq0}) has no semi-singular solutions, if N=4
$\bb\ge3\max\{\mu_1, \mu_2\}$. Then for any a solution $(u, v)$ of (\ref{eq0}) with $K(u, v)<0$, there exist constants $\mathcal{C}_2>\mathcal{C}_1>0$ such that
\be\label{eq8}\mathcal{C}_1|x|^{\frac{2-N}{2}}\le u(x), v(x)\le \mathcal{C}_2|x|^{\frac{2-N}{2}},\quad\forall\,x\in\RN\setminus\{0\}.\ee
\end{itemize}\et

In particular, Theorems \ref{th1}-\ref{th2} and Corollary \ref{corollary1-1} give a classification of solutions to system (\ref{eq0}) in the case $N\ge 5$. Remark that, the above results are very important when we study system (\ref{eq0}) in a {\it punctured ball}. We will treat this problem in a subsequent work \cite{ChenLin2} where, among other things, the {\it asymptotic symmetry} of singular solutions will be studied. Remark that the scalar equation (\ref{eq4}) in a punctured ball was well investigated in the cerebrated paper \cite{CGS}; see also \cite{CL1, CL2, CL3, KMPS, Lin, TZ} for subsequent progress. Our study here and the subsequent \cite{ChenLin2} can be seen as extensions of \cite{CGS, Fowler} to systems.

Theorem \ref{th2} indicates that, for $N\ge 4$, $0$ is a non-removable singularity of $u$ if and only if $0$ is a non-removable singularity of $v$. So far, we have no idea whether this conclusion holds or not for $N=3$, but we tend to believe that it is also true for $N=3$,
by which we can obtain the sharp global estimates of both-singular solutions.
The case $N=4$ is quite different, and we can not prove the sharp global estimates of both-singular solutions via the nonexistence of semi-singular solutions. Remark that, the sharp global estimates of singular solutions for the scalar equation (\ref{eq4}) is natural, and the proof is very simple. However, our problem (\ref{eq0}) is a coupled system, which turns out to be much more difficult and complicated than (\ref{eq4}). Moreover, in the equation of $u$ in (\ref{eq0}), the power of $u$ in the coupling term is $2^\ast/2-1=\frac{2}{N-2}$, which $>1$ if $N=3$, $=1$ if $N=4$ and $<1$ if $N\ge5$. Clearly this fact makes the study of (\ref{eq0}) depending heavily on the dimensions. This is a further difference between (\ref{eq0}) and (\ref{eq4}). Here
we need a technical assumption $\bb\ge 3\max\{\mu_1, \mu_2\}$ for $N=4$.  We conjecture that (\ref{eq8}) holds for all both-singular solutions not only in the case $N=4$ with $0<\bb<3\max\{\mu_1, \mu_2\}$ but also in the case $N=3$ without assuming the nonexistence of semi-singular solutions.

Let $V$ be a singular solution of (\ref{eq4}) with $\lim_{|x|\to 0}V(x)=+\iy$, then $(kV, lV)$ is a both-singular solution of (\ref{eq0}), where $k, l>0$ satisfies (\ref{eq3}). As pointed out before, there is a basic question remaining: {\it whether any both-singular solutions are of the form $(kV, lV)$ just as entire solutions in Theorem A}? This question seems very tough and remains completely open. Remark that, although we will prove the radial symmetry of a singular solution $(u, v)$ (either semi-singular or both-singular) in Theorem \ref{th2-1} below, we can not obtain the conclusion $u/v\equiv\text{constant}$ via the ideas of proving Theorem A in \cite{ChenLi, GL}. The essential difference is that system (\ref{eq1}) is invariant under translations (this property, together with Kelvin transform, plays the key role in Theorem A)
but our problem (\ref{eq0}) is {\it not}. Therefore, system (\ref{eq0}) can not be reduced to a single equation. In fact, we will prove our results by analyzing an ODE system, which turns out to be very delicate.

The rest of this paper is organized as follows. In Section 2, we prove that singular solutions of (\ref{eq0}) are radially symmetric via the method of moving planes. We mainly follow the idea in \cite{CGS} where the radial symmetry of singular solutions of (\ref{eq4}) was proved. In Section 3, we establish some crucial lemmas by studying the corresponding ODE system. Theorems \ref{th1} and \ref{th2} are proved in Sections 4 and 5 respectively. As pointed out before, since our problem is a coupled system, we will see that the arguments are completely different and much more delicate comparing to the scalar case, and some new ideas are needed. In the sequel, we denote positive constants (possibly different in different places) by $C, C_0, C_1,\cdots$.

\section{Radial symmetry}
\renewcommand{\theequation}{2.\arabic{equation}}

In this section, we use the method of moving planes to prove the radial symmetry of positive singular solutions, which is an extension of a classical result from Caffarelli, Gidas and Spruck \cite{CGS}.
Remark that, since we will use the Kelvin transform, it does not seem trivial to obtain the strictly decreasing property of $u$ and $v$ as functions of $r=|x|$ directly from this argument. Therefore, instead of giving the proof here, we will prove this strictly decreasing property via a different observation in Corollary \ref{corollary3-1} in Section 3.

\bt\label{th2-1} Let $(u, v)$ be a singular solution (either both-singular or semi-singular) of (\ref{eq0}), then both $u$ and $v$ are radially symmetric about the origin and are strictly decreasing with respect to $r=|x|>0$.\et

\bp
Without loss of generality, we assume that $0$ is a non-removable singularity of $u$, namely $\lim_{|x|\to 0}u(x)=+\iy$.
We mainly follow the argument in \cite{CGS} to apply the method of moving planes.
Fix an arbitrary point $z\neq 0$. Define the Kelvin transform
$$U(x)=|x|^{2-N}u\left(z+\frac{x}{|x|^2}\right),\quad V(x)=|x|^{2-N}v\left(z+\frac{x}{|x|^2}\right).$$
Then $U, V$ are singular at $0$ and $z_0=-z/|z|^2$, and
\be\label{eq2-1}
\begin{cases}
\begin{split}-\Delta U =
\mu_1 U^{2^\ast-1}+\beta U^{\frac{2^\ast}{2}-1} V^{\frac{2^\ast}{2}}\\
-\Delta U =\mu_2 V^{2^\ast-1}+\beta V^{\frac{2^\ast}{2}-1}U^{\frac{2^\ast}{2}}\\
\end{split}  \quad   \text{in $\RN\setminus\{0, z_0\}$},\\
U, V> 0 \quad\hbox{and}\quad
U, V\in C^2(\RN\setminus\{0, z_0\}).\end{cases}\ee
We shall prove that both $U$ and $V$ are axisymmetric with respect to the axis going through $0$ and $z$.
For this purpose, we consider any reflection direction $\tau$ orthogonal to this axis. Without loss of generality, we may assume that $\tau=(0,\cdots, 0, 1)$ is the positive $x_N$ direction. For $\la>0$ we consider the reflection
$$x=(x_1, \cdots, x_{N-1}, x_N)\mapsto x^\la=(x_1,\cdots,x_{N-1}, 2\la-x_N),$$
where $x\in \Sg_\la:=\{x\in\RN : x_N>\la\}$. Since $U, V$ have the harmonic asymptotic expansion \cite[(2.6)]{CGS} at $\iy$, by \cite[Lemma 2.3]{CGS}, there exist large positive constants $\bar{\la}>10$ and  $R>|z_0|+10$ such that for any $\la\ge \bar{\la}$, we have
\be\label{eq2-2}U(x)<U(x^\la),\quad V(x)<V(x^\la),\quad\text{for $x\in \Sg_\la$ and $|x^\la|>R$}.\ee
Since the maximum principle holds for super harmonic functions with isolated singularities (cf. \cite[Lemma 2.1]{CL1}), there exists $C>0$ such that
\be\label{eq2-3}U(x), V(x)\ge C,\quad\text{for}\,\,\, x\in \overline{B_R}\setminus\{0, z_0\}.\ee
Since $U(x), V(x)\to 0$ as $|x|\to +\iy$, by (\ref{eq2-2})-(\ref{eq2-3}) there exists $\la_0>\bar{\la}$ such that
for any $\la\ge \la_0$, we have
\be\label{eq2-4}U(x)\le U(x^\la),\quad V(x)\le V(x^\la),\quad\text{for $x\in \Sg_\la$ and $x^\la\not\in \{0, z_0\}$}.\ee
Let $\la^\ast:=\inf\{\tilde{\la}>0 \,\,|\,\, (\ref{eq2-4})\,\, \text{holds for all $\la\ge \tilde{\la}$}\}$. Suppose $\la^\ast>0$. Then (\ref{eq2-4}) holds for $\la=\la^*$. Since $\lim_{|x|\to 0}u(x)=+\iy$, we see that $\lim_{x\to z_0}U(x)=+\iy$ and so $U(x)\not\equiv U(x^{\la^*})$. Then by applying the maximum principle to $U(x^{\la^*})-U(x)$ and $V(x^{\la^*})-V(x)$ in (\ref{eq2-1}), we easily conclude that
\be\label{eq2-5}U(x)< U(x^{\la^\ast}),\quad V(x)< V(x^{\la^*}),\quad\text{for $x\in \Sg_{\la^*}$ and $x^{\la^*}\not\in \{0, z_0\}$}.\ee
This together with the Hopf boundary lemma yield for $x\in\partial\Sg_{\la^*}$ (note that $0, z_0\not\in \partial\Sg_{\la^*}$ because of $\la^\ast>0$) that
\be\label{eq2-6}\frac{\partial(U(x^{\la^\ast})-U(x))}{\partial x_N}=-2\frac{\partial U(x)}{\partial x_N}>0,\,\,\, \frac{\partial(V(x^{\la^\ast})-V(x))}{\partial x_N}=-2\frac{\partial V(x)}{\partial x_N}>0.\ee
By the definition of $\la^*$, there exists $\la_j\uparrow \la^*$ such that (\ref{eq2-4}) does not hold for $\la=\la_j$. Without loss of generality and up to a subsequence, we may assume the existence of $x_j\in \Sg_{\la_j}$ such that $U(x_j^{\la_j})<U(x_j)$. By \cite[Lemma 2.4]{CGS} (the plane $x_N=0$ there corresponds to $x_N=\la^*$ here), we know that $|x_j|$ are uniformly bounded. Up to a subsequence, $x_j\to \bar{x}\in\overline{\Sg_{\la^\ast}}$ with $U(\bar{x}^{\la^\ast})\le U(\bar{x})$. By (\ref{eq2-5}) we know $\bar{x}\in\partial\Sg_{\la^*}$ and then $\frac{\partial U}{\partial x_N}(\bar{x})\ge 0$, a contradiction with (\ref{eq2-6}). Therefore, $\la^*=0$ and so both $U$ and $V$ are axisymmetric with respect to the axis going through $0$ and $z$. Since $z$ is arbitrary, $u, v$ are both radially symmetric about the origin.
The strictly decreasing property of $u$ and $v$ as functions of $r=|x|$ will be proved via a different observation in Corollary \ref{corollary3-1} in Section 3.\ep

\section{Useful observations via ODE analysis}
\renewcommand{\theequation}{3.\arabic{equation}}

In this section, we establish some crucial lemmas via ODE analysis. Let $(u, v)$ be a solution of (\ref{eq0}) with $K(u, v)\in\R$. By Theorem A and Theorem \ref{th2-1}, we may assume that $u(x)=u(|x|)$ and $v(x)=v(|x|)$ are radially symmetric functions. Denote $\dd=\frac{N-2}{2}$ and $p=\frac{2^\ast}{2}$ for convenience. Now we use the classical change of variables from Fowler \cite{Fowler}. Let $t=-\ln r$ and
$$w_1(t):=r^{\dd}u(r)=e^{-\dd t}u(e^{-t}),\quad w_2(t):=r^{\dd}v(r)=e^{-\dd t}v(e^{-t}).$$
Then it is easy to see from (\ref{eq0}) that $(w_1, w_2)$ satisfies
\be\label{eq3-1}
\begin{cases}
\begin{split}w_{1}''-\dd^2w_{1}+\mu_1w_1^{2p-1}+\bb w_1^{p-1}w_2^{p}=0\\
w_{2}''-\dd^2w_{2}+\mu_2w_2^{2p-1}+\bb w_2^{p-1}w_1^{p}=0\\
\end{split}  \quad   t\in\R,\\
w_1, w_2> 0 \quad\hbox{and}\quad
w_1, w_2\in C^2(\R),\end{cases}\ee
where $w_{i}':=\frac{d w_i}{dt}$ and $w_{i}'':=\frac{d^2 w_i}{dt^2}$. Define
{\allowdisplaybreaks
\begin{align}\label{eq3-2}
\Psi(t):=&\frac{1}{2}(|w_1'|^2+|w_2'|^2-\dd^2w_1^2-\dd^2w_2^2)(t)\nonumber\\
&+\frac{1}{2p}(\mu_1|w_1|^{2p}+2\bb |w_1|^p|w_2|^p+\mu_2|w_2|^{2p})(t).
\end{align}
}%
Then $\Psi'(t)\equiv 0$, namely $\Psi(t)$ is a constant. In fact, by a direct computation, we easily deduce from
(\ref{eq7}) that $K(r; u, v)=\sg_{N-1}\Psi(t)$, where $t=-\ln r$ and $\sg_{N-1}$ is the area of the unit sphere in $\RN$. Hence
\be\label{eq3-3}\Psi(t)\equiv \frac{K(u, v)}{\sg_{N-1}}=:K.\ee
Combining this with (\ref{eq3-2}), it is easy to prove the existence of $C>0$ such that $w_1, w_2, |w_{1}'|, |w_{2}'|\le C$ for all $t\in \R$.

For the scalar curvature equation (\ref{eq4}), the corresponding constant-value function is $\Phi(t)=\frac{1}{2}|w'|^2-\frac{1}{2}\dd^2w^2+\frac{1}{2^\ast}w^{2^\ast}$, where $w(t)=e^{-\dd t}u(e^{-t})$ and $u$ is a solution of (\ref{eq4}) (we may assume that $u$ is radially symmetric by the classical result of Caffarelli, Gidas and Spruck \cite{CGS}). The qualitative analysis of $\Phi(t)$ is easy, from which one can describe all solutions of (\ref{eq4}); see \cite{Fowler} or the Introduction of \cite{CGS} for example. However, since our problem (\ref{eq0}) is a system, the
qualitative analysis of $\Psi(t)$ is completely different
and much more delicate comparing to that of $\Phi(t)$.

\bl\label{lemma3-1}Fix $i\in \{1, 2\}$. If $w_i$ has a limit $C_i$ at $+\iy$ or $-\iy$, then $C_i\le (\dd^2/\mu_i)^{\frac{1}{2p-2}}$. Moreover, $C_i=0$ provided $K\ge 0$.\el

\bp Without loss of generality, we assume $\lim_{t\to+\iy}w_1=C_1$. Suppose $C_1>0$. Up to a subsequence, by (\ref{eq3-1}) we may assume that
$w_1(\cdot+n)\to C_1$ and $w_2(\cdot+n)\to v_2\ge 0$ uniformly in $C^2_{loc}(\R)$, and
$-\dd^2 C_1+\mu_1C_1^{2p-1}+\bb C_1^{p-1}v_2^{p}=0$. This implies $C_1\le (\dd^2/\mu_1)^{\frac{1}{2p-2}}$. Moreover, $v_2\equiv C_2$ is also a constant and
{\allowdisplaybreaks
\begin{align*}
&-\dd^2 C_2+\mu_2C_2^{2p-1}+\bb C_2^{p-1}C_1^{p}=0,\\
&-\frac{1}{2}(\dd^2C_1^2+\dd^2C_2^2)+\frac{1}{2p}(\mu_1C_1^{2p}+2\bb C_1^pC_2^p+\mu_2C_2^{2p})=K.
\end{align*}
}%
These imply $K=-\frac{1}{N}(\dd^2C_1^2+\dd^2C_2^2)<0$. Hence $C_1=0$ provided $K\ge 0$.
\ep

Define two functions $f_1, f_2: \R\to \R$ by
\be\label{eq3-4}f_1(t):=-\frac{1}{2}|w_{1}'|^2+\frac{\dd^2}{2}w_1^2-\frac{\mu_1}{2p}w_1^{2p},\,\,\, f_2(t):=-\frac{1}{2}|w_{2}'|^2+\frac{\dd^2}{2}w_2^2-\frac{\mu_2}{2p}w_2^{2p}.\ee
By (\ref{eq3-1}) we have $f_1'(t)=\bb w_1^{p-1}w_2^p w_{1}'$ and $f_2'(t)=\bb w_2^{p-1}w_1^p w_{2}'$. That is, for $i=1, 2$,
{\it the monotonicity of $f_i$ is exactly the same as the monotonicity of $w_i$}. This monotonicity property is very important and will be used frequently in the sequel. The following two lemmas are our important observations via this monotonicity property, and will play crucial roles in the proof of Theorem \ref{th1}.
\bl\label{lemma3-2}
$$w_i(t)<\la_i:=\left(\frac{p\dd^2}{\mu_i}\right)^{\frac{1}{2p-2}},\quad \forall\, t\in\R\,\,\,\text{and}\,\,\, i=1, 2.$$ \el

\bp
Suppose by contradiction that $w_1(s)\ge \la_1> (\dd^2/\mu_1)^{\frac{1}{2p-2}}$ for some $s\in\R$. By Lemma \ref{lemma3-1} there exists a local maximum point $t_0\in \R$ of $w_1$ such that $w_1(t_0)\ge w_1(s)\ge\la_1$. Then
\be\label{eq3-5}f_1(t_0)=\frac{\dd^2}{2}w_1^2(t_0)-\frac{\mu_1}{2p}w_1^{2p}(t_0)\le 0.\ee
If $w_1'(t)<0$ for all $t>t_0$, it follows from Lemma \ref{lemma3-1} that $w_1\downarrow C_1\le (\dd^2/\mu_1)^{\frac{1}{2p-2}}$ and $f_1(t)\downarrow- C<0$ as $t\uparrow+\iy$. By (\ref{eq3-4}), $w_1'(t)\to 0$ as $t\uparrow +\iy$ and so
$$0>-C=\frac{\dd^2}{2}C_1^2-\frac{\mu_1}{2p}C_1^{2p}\ge 0,$$
a contradiction. Hence there exists $t_1>t_0$ such that $w_1'(t)<0$ for $t\in (t_0, t_1)$ and $w_1'(t_1)=0$. Then $w_1''(t_1)\ge 0$ and
\be\label{eq3-6}f_1(t_1)=\frac{\dd^2}{2}w_1^2(t_1)-\frac{\mu_1}{2p}w_1^{2p}(t_1)<f_1(t_0)\le 0,\ee
namely $w_1(t_1)>\la_1$. But this implies from (\ref{eq3-1}) that $w_1''(t_1)<0$, a contradiction.
\ep

\bl\label{lemma3-3} $f_i(t)>0$ and so $|w_i'(t)|<\dd w_i(t)$ for all $t\in\R$ and $i=1, 2$. \el

\bp
Suppose $f_1(s)\le 0$ for some $s\in\R$. If $f_1'(s)=0$, then $w_1'(s)=0$ and so
(\ref{eq3-5}) holds for $s$, which implies $w_1(s)\ge \la_1$, a contradiction with Lemma \ref{lemma3-2}.
So we may assume, without loss of generality, that $f_1'(s)>0$, namely $w_1'(s)>0$. If $w_1'(t)>0$ for all $t\le s$, then $w_1(t)\downarrow C_1\le (\dd^2/\mu_1)^{\frac{1}{2p-2}}$ and $f_1(t)\downarrow -C<0$ as $t\downarrow-\iy$, and we get a contradiction as in the proof of Lemma \ref{lemma3-2}. So there exists $t_0<s$ such that $w_1'(t)>0$ for $t\in (t_0, s]$ and $w_1'(t_0)=0$. Similarly as (\ref{eq3-6}) we get $w_1(t_0)>\la_1$, a contradiction with Lemma \ref{lemma3-2}.
\ep

As an application of Lemma \ref{lemma3-3}, we finish the proof of Theorem \ref{th2-1} by proving the following corollary.

\begin{corollary}\label{corollary3-1} Both $u'(r)<0$ and $v'(r)<0$ for all $r>0$.  \end{corollary}

\bp Since $u(r)=r^{-\dd}w_1(-\ln r)$, Lemma \ref{lemma3-3} yields $u'(r)=-r^{-\dd-1}(w_1'(t)+\dd w_1(t))<0$, where $t=-\ln r$.
\ep

\section{Removable singularity}
\renewcommand{\theequation}{4.\arabic{equation}}

In this section, we always assume $K(u, v)\ge 0$ (except in Corollary \ref{corollary4-1}). We want to prove that $K(u, v)=0$ and $u, v\in C^2(\RN)$. Obviously, Theorem \ref{th1} follows immediately. First we consider the simple case $N\le 4$.

\begin{proof}[Proof of Theorem \ref{th1} for $N\le 4$] Recall $K\ge 0$. It follows from (\ref{eq3-2}), (\ref{eq3-4}) and Lemma \ref{lemma3-3} that
\be\label{eq3-7}0<f_1(t)+ f_2(t)=\frac{\bb}{p}w_1^pw_2^p-K\le \frac{\bb}{p}w_1^pw_2^p.\ee
Lemma \ref{lemma3-1} implies the existence of $s\in \R$ such that $w_1'(s)=0$.  Let $s\in \R$ be any a point such that $w_1'(s)=0$. Then Lemma \ref{lemma3-3} and (\ref{eq3-7}) give
$$f_1(s)=\frac{\dd^2}{2}w_1^2(s)-\frac{\mu_1}{2p}w_1^{2p}(s)<\frac{\bb}{p}w_1^p(s)w_2^p(s),$$
and so
{\allowdisplaybreaks
\begin{align*}
w_1''(s)&=\dd^2w_1(s)-\mu_1w_1^{2p-1}(s)-\bb w_1^{p-1}(s)w_2^{p}(s)\\
&<\frac{1-p}{p}\mu_1w_1^{2p-1}(s)+\frac{2-p}{p}\bb w_1^{p-1}(s)w_2^{p}(s).
\end{align*}
}%
Now we assume $N\le 4$. Then $p=\frac{2^\ast}{2}\ge 2$ and so $w_1''(s)<0$, that is, $s$ must be a local maximum point of $w_1$, and $w_1$ has no local minimum points. Therefore, $w_1$ has a unique maximum point $t_1$, $w_1'(t)>0$ for $t<t_1$ and $w_1'(t)<0$ for $t>t_1$. Similarly, $w_2$ has a unique maximum point $t_2$, $w_2'(t)>0$ for $t<t_2$ and $w_2'(t)<0$ for $t>t_2$. Lemma \ref{lemma3-1} gives that $w_1\downarrow 0$ and $w_2\downarrow 0$ as $t\uparrow +\iy$. This together with Lemma \ref{lemma3-3} yields $K=0$, namely $K(u, v)=0$. Let $w=w_1+w_2$ and
$y(t)=w'(t)+\dd w(t)$, then $w'(t)<0$ for all $t>\max\{t_1, t_2\}$ and
$$y'-\dd y=w''-\dd^2 w\ge -C w^{2p-1}.$$
Since $y>0$ is bounded, it is easy to prove that
$$w'(t)+\dd w(t)=y(t)\le\frac{C}{\dd}w^{2p-1}(t),\quad\forall\,t> \max\{t_1, t_2\},$$
and then $(e^{\dd t}w(t))^{2-2p}-\frac{C}{\dd^2}e^{(2-2p)\dd t}$ is strictly increasing for $t> \max\{t_1, t_2\}$. Combining this with $w(t)\downarrow 0$ as $t\uparrow+\iy$, we get that $e^{\dd t}w(t)\le C$ uniformly for $t>0$ large. That is, $u(r)+v(r)\le C$ uniformly for $r>0$ small. Therefore, $u, v\in C^2(\RN)$ from the standard elliptic regularity theory. \end{proof}

In the rest of this section, we consider the case $N\ge 5$, where $1<p<2$ and so the above proof does not apply.
Therefore we need to develop completely different ideas.  The following proof is very delicate and long, and we divide it by proving some lemmas.

\bl\label{lemma3-4} If $t_0$ be a local minimum point of $w_1$, then
\be\label{eq3-8}w_1(t_0)<\la_1^*:=\left(\frac{\dd^2}{\mu_1}\right)^{\frac{1}{2p-2}},\quad f_2(t_0)<\frac{N-4}{N}f_1(t_0).\ee
Similarly, if $t_0$ is a local minimum point of $w_2$, then
$$w_2(t_0)<\la_2^*:=\left(\frac{\dd^2}{\mu_2}\right)^{\frac{1}{2p-2}},\quad f_1(t_0)<\frac{N-4}{N}f_2(t_0).$$
 \el

\bp We give the proof of (\ref{eq3-8}). The inequality $w_1(t_0)<\la_1^*$ follows directly from $w_1''(t_0)\ge 0$.
Suppose $f_2(t_0)\ge\frac{N-4}{N}f_1(t_0)$, then (\ref{eq3-7}) implies
$$\frac{\bb}{p}(w_1^pw_2^p)(t_0)-K=f_1(t_0)+ f_2(t_0)\ge\frac{2N-4}{N}f_1(t_0),$$
and so $\dd^2 w_1(t_0)\le \frac{1}{p}\mu_1w_1^{2p-1}(t_0)+\bb (w_1^{p-1}w_2^p)(t_0)$, which is a contradiction with $w_1''(t_0)\ge 0$.\ep

\bl\label{lemma3-5} Both $\liminf\limits_{t\to\pm\iy } w_1=0$ and $\liminf\limits_{t\to\pm\iy}w_2(t)=0$.\el

\bp
Assume by contradiction that $\liminf\limits_{t\to+\iy } w_1=C_1>0$. If there exists $t_n\uparrow+\iy$ such that $w_2(t_n)\to 0$, then, up to a subsequence,
$w_1(\cdot+t_n)\to v_1\ge C_1$ and $w_2(\cdot+t_n)\to v_2\ge 0$ uniformly in $C^2_{loc}(\R)$, and
\be\label{eq3-9}
\begin{cases}
\begin{split}v_{1}''-\dd^2v_{1}+\mu_1v_1^{2p-1}+\bb v_1^{p-1}v_2^{p}=0\\
v_{2}''-\dd^2v_{2}+\mu_2v_2^{2p-1}+\bb v_2^{p-1}v_1^{p}=0\\
\end{split}  \quad   t\in\R.\end{cases}\ee
Since $v_2(0)=0$, the maximum principle yields $v_2\equiv 0$, that is,
\be\label{eq4-1}v_{1}''-\dd^2v_{1}+\mu_1v_1^{2^\ast-1}=0 \quad\text{and}\quad v_1\ge C_1.\ee
Then by (\ref{eq3-2}) and a classical result from Fowler \cite{Fowler} (i.e. Theorem B), we get
\be\label{eq4-2}K\equiv \frac{1}{2}|v_1'(t)|^2-\frac{1}{2}\dd^2v_1^2(t)+\frac{1}{2^\ast}\mu_1v_1^{2^\ast}(t)<0,\ee
a contradiction with $K\ge 0$. Hence $\liminf\limits_{t\to+\iy } w_2=C_2>0$.

By Taking $t_n\to+\iy$ such that $w_1(t_n)\to C_1$ and repeating the progress of (\ref{eq3-9}), without loss of generality, we may assume in the beginning that
\be\label{eq3-10}w_1(0)=\inf_{t\in\R}w_1(t)=C_1 \quad\text{and}\quad \inf_{t\in\R}w_2(t)= C_2.\ee
We claim that
\be\label{eq3-11}f_1(0)=\inf_{t\in\R}f_1(t).\ee

Take any $t\in\R$. By Lemma \ref{lemma3-1} and the same monotonicity of $w_1$ and $f_1$, there exists a local minimum point $t_1$ of $w_1$ such that $w_1(t)\ge w_1(t_1)$ and $f_1(t)\ge f_1(t_1)$. Lemma \ref{lemma3-4} tells us that $C_1=w_1(0)\le w_1(t_1)<\la_1^\ast$. Since the function $\vp(s)=\frac{1}{2}\dd^2s^2-\frac{1}{2p}\mu_1s^{2p}$ is strictly increasing for $s\in [0, \la_1^\ast]$, we have
\be\label{eq3-13}f_1(0)=\vp(w_1(0))\le\vp(w_1(t_1))=f_1(t_1)\le f_1(t),\ee
and so (\ref{eq3-11}) holds. By Lemma \ref{lemma3-4} we know
\be\label{eq3-12}\inf_{t\in\R}f_2(t)\le f_2(0)<\frac{N-4}{N}f_1(0)=\frac{N-4}{N}\inf_{t\in\R}f_1(t).\ee
If there exists $t_2\in\R$ such that $w_2(t_2)=\inf w_2=C_2$, by repeating the above argument we see that $\inf_{t\in\R}f_1(t)\le f_1(t_2)<\frac{N-4}{N}f_2(t_2)=\frac{N-4}{N}\inf_{t\in\R}f_2(t)$, a contradiction.
Hence $C_2=\inf w_2$ can not be attained. By Lemma \ref{lemma3-1} there exists a sequence of local minimum points $t_n$ of $w_2$ such that $|t_n|\to\iy$ and $w_2(t_n)\to C_2$.
Up to a subsequence,
$w_1(\cdot+t_n)\to v_1\ge C_1$ and $w_2(\cdot+t_n)\to v_2\ge C_2$ uniformly in $C^2_{loc}(\R)$, where $(v_1, v_2)$ satisfies (\ref{eq3-9}) and $v_2(0)=C_2=\inf_{t\in\R}v_2(t)$.
As before, we define $g_1, g_2: \R\to \R$ by
$$g_1(t):=-\frac{1}{2}|v_{1}'|^2+\frac{\dd^2}{2}v_1^2-\frac{\mu_1}{2p}v_1^{2p},\,\,\, g_2(t):=-\frac{1}{2}|v_{2}'|^2+\frac{\dd^2}{2}v_2^2-\frac{\mu_2}{2p}v_2^{2p}.$$
Clearly Lemma \ref{lemma3-4} also holds for $v_i$ and $g_i$. Hence
$$g_2(0)=\frac{\dd^2}{2}C_2^2-\frac{1}{2p}\mu_2 C_2^{2p}>\frac{N}{N-4}g_1(0).$$
Since $\inf w_2= C_2$, the same argument as (\ref{eq3-13}) yields $\frac{\dd^2}{2}C_2^2-\frac{1}{2p}\mu_2 C_2^{2p}\le\inf_{t\in\R}f_2(t)$, namely $g_2(0)\le \inf_{t\in\R} f_2(t)$. Therefore, we deduce from (\ref{eq3-12}) that
\begin{align*}
g_2(0)< \frac{N-4}{N}\inf_{t\in\R}f_1(t)\le\frac{N-4}{N}\lim_{t_n\to \iy}f_1(t_n)=\frac{N-4}{N}g_1(0),
\end{align*}
a contradiction. This completes the proof.
\ep

A direct corollary of Lemma \ref{lemma3-5} is the following result Corollary \ref{corollary4-1}. Remark that we need assumption $K\ge 0$ in the proof of Lemma \ref{lemma3-5}, but
Corollary \ref{corollary4-1} holds for all $K\in \R$ provided
that $N\ge 5$ and Lemma \ref{lemma3-5} hold. We will also use Corollary \ref{corollary4-1} in
the proof of Theorem \ref{th2} where $K<0$.

\begin{corollary}\label{corollary4-1} Assume that $N\ge 5$ and Lemma \ref{lemma3-5} hold. Then
\be\label{eq4-3}\liminf_{t\to\pm \iy}(w_1+w_2)(t)=0\quad\text{and}\quad K(u, v)=0.\ee
\end{corollary}

\bp
If there exists $T$ such that $w_1'(t)<0$ for all $t>T$, then Lemma \ref{lemma3-5} (not by Lemma \ref{lemma3-1} since we don't assume $K\ge 0$ here) gives that $w_1(t)\to 0$ as $t\to+\iy$, and so $\liminf_{t\to+\iy}(w_1+w_2)(t)=0$. Otherwise, there exists a sequence of local minimum point $t_n$ of $w_1$ such that $t_n\to+\iy$ and $w_1(t_n)\to 0$. By $w_1''(t_n)\ge 0$ we obtain
$\bb w_2^p(t_n)<\dd^2 w_1^{2-p}(t_n)$. Since $1<p<2$ from $N\ge 5$, we have $w_2(t_n)\to 0$ and so $(w_1+w_2)(t_n)\to 0$. Therefore the first identity of (\ref{eq4-3}) holds.

To prove $K=0$, we take $t_n\to+\iy$ such that $(w_1+w_2)(t_n)\to 0$. Up to a subsequence, $w_1(\cdot+t_n)\to v_1\ge 0$ and $w_2(\cdot+t_n)\to v_2\ge 0$ uniformly in $C^2_{loc}(\R)$, where $v_1(0)=v_2(0)=0$ and $(v_1, v_2)$ satisfies (\ref{eq3-9}). Then the maximum principle yields $v_1\equiv v_2\equiv 0$ and so $K=0$.\ep

Now we improve Lemma \ref{lemma3-5} by proving the following lemma.

\bl\label{lemma4-1}There holds $\lim_{t\to+\iy}\min\{w_1(t), w_2(t)\}=0$. Moreover, if both $w_1$ and $w_2$ have a sequence of local minimum points which tends to $+\iy$, then there exists an increasing sequence $T_n (n\ge 1)$ such that $T_n\to +\iy$ and for any $k\ge 1$,
\begin{itemize}
\item[$(1)$] $T_{2k-1}$ is a local minimum point of $w_1$ and $T_{2k}$ is a local minimum point of $w_2$;
\item[$(2)$] $w_1'(t)<0$ on $[T_{2k}, T_{2k+1})$ and $w_2'(t)<0$ on $[T_{2k-1}, T_{2k})$.
\end{itemize}
\el

\bp
If there exists $T>0$ such that $w_1'(t)<0$ for all $t>T$, we are done since Lemma \ref{lemma3-1} yields that $w_1(t)\to 0$ as $t\to +\iy$. Hence, we may assume that both $w_1$ and $w_2$ have a sequence of local minimum points which tends to $+\iy$.  Recall that the monotonicity of $f_i$ is exactly the same as that of $w_i$, and we will use this fact frequently in the following proof.

Assume by contradiction that
\begin{itemize}
\item[$(H)$] $w_1'(t)>0$ whenever $t$ is a local minimum point of $w_2$, and  $w_2'(t)>0$ whenever $t$ is a local minimum point of $w_1$.
\end{itemize}

Let $t_1$ be a local minimum point of $w_1$, then $w_2'(t_1)>0$ by assumption $(H)$ and Lemma \ref{lemma3-4} gives
\be\label{eq4-5}f_1(t_1)>\frac{N}{N-4}f_2(t_1).\ee
Clearly there exist $t_2>\tau_1>t_1$ such that $w_2'(t)>0$ in $[t_1, \tau_1)$, $w_2'(t)<0$ in $(\tau_1, t_2)$, $\tau_1$ is a local maximum point of $w_2$ and $t_2$ be a local minimum point of $w_2$. By assumption $(H)$ we have $w_1'(t_2)>0$. Hence there exists $\rho_1\in [t_1, t_2)$ such that $w_1'(t)>0$ in $(\rho_1, t_2]$ and $\rho_1$ is a local minimum point of $w_1$ (maybe $\rho_1=t_1$). Then Lemma \ref{lemma3-4} gives
$$f_2(t_2)>\frac{N}{N-4}f_1(t_2)\quad\text{and}\quad f_1(\rho_1)>\frac{N}{N-4}f_2(\rho_1).$$
Then there exists $\rho\in (\rho_1, t_2)$ such that $f_2(\rho)=f_1(\rho)$ and $f_2(t)>f_1(t)$ for all $t\in (\rho, t_2]$. In particular, $f_2'(\rho)\ge f_1'(\rho)>0$ (since $w_1'(t)>0$ on $(\rho_1, t_2]$), namely $w_2'(\rho)>0$ and so $\rho_1<\rho<\tau_1$. In a word, $t_1\le\rho_1<\tau_1<t_2$. Therefore, by Lemma \ref{lemma3-4} and the monotonicity of $f_i$, we have
{\allowdisplaybreaks
\begin{align}\label{eq4-61}f_2(t_2)&>\frac{N}{N-4}f_1(t_2)>\frac{N}{N-4}f_1(\rho_1)\nonumber\\
&>\left(\frac{N}{N-4}\right)^2f_2(\rho_1)\ge\left(\frac{N}{N-4}\right)^2f_2(t_1).\end{align}
}%
Since $t_2$ is a local minimum point of $w_2$ and $w_1'(t_2)>0$, by repeating the above argument from (\ref{eq4-5}) to (\ref{eq4-61}) (change the role of $w_1$ and $w_2$), we can prove the existence of $t_3>t_2$ such that $t_3$ is a local minimum point of $w_1$ and $w_2'(t_3)>0$ and
\begin{align}\label{eq4-7}f_1(t_3)>\left(\frac{N}{N-4}\right)^2f_1(t_2)\quad\text{(similarly as (\ref{eq4-61}))}.\end{align}
This together with (\ref{eq4-61}) give
\be\label{eq4-6}f_1(t_3)>\left(\frac{N}{N-4}\right)^3f_2(t_1).\ee
From (\ref{eq4-5}), (\ref{eq4-61}) and (\ref{eq4-7}), we can continue this progress to obtain an increasing sequence $t_n$ with
$$f_1(t_{2n+1})\ge \left(\frac{N}{N-4}\right)^{2n+1}f_2(t_1)\quad\text{and}\quad f_2(t_{2n})\ge \left(\frac{N}{N-4}\right)^{2n}f_2(t_1),$$
which contradicts the fact that $f_1$ and $f_2$ are both uniformly bounded in $\R$. Hence assumption $(H)$ is not true.

Without loss of generality, we may assume the existence of $t_0$ such that $t_0$ is a local minimum point of $w_1$ and $w_2'(t_0)\le 0$. Then Lemma \ref{lemma3-4} gives
\be\label{eq4-8}f_2(t_0)<\frac{N-4}{N}f_1(t_0).\ee
This also means that $t_0$ is not a local minimum point of $w_2$. Consequently there exists $T_0>t_0$ such that
$T_0$ is a local minimum point of $w_2$ and $w_2'(t)<0$ on $(t_0, T_0)$. Again
\be\label{eq4-9}f_1(T_0)<\frac{N-4}{N}f_2(T_0).\ee
By (\ref{eq4-8})-(\ref{eq4-9}), there exists $\nu_0\in (t_0, T_0)$ such that $f_2(\nu_0)=f_1(\nu_0)$ and $f_2(t)>f_1(t)$ for $t\in (\nu_0, T_0]$. Then $f_1'(\nu_0)\le f_2'(\nu_0)<0$ (since $w_2'(\nu_0)<0$), namely $w_1'(\nu_0)<0$. By $f_2(t)>f_1(t)$ for $t\in (\nu_0, T_0]$, we see that $w_1$ has no local minimum points in $[\nu_0, T_0]$ and so $w_1'(T_0)<0$. Moreover,
\be\label{eq4-10}f_2(T_0)<f_2(t_0)<\frac{N-4}{N}f_1(t_0).\ee
Since $T_0$ is a local minimum point of $w_2$ and $w_1'(T_0)<0$, by repeating the above argument from (\ref{eq4-8}) to (\ref{eq4-10}) (change the role of $w_1$ and $w_2$), we can prove the existence of $T_1>T_0$ such that $T_1$ is a local minimum point of $w_1$, $w_1'(t)<0$ on $[T_0, T_1)$, $w_2'(T_1)<0$ and
\begin{align}\label{eq4-11}f_1(T_1)<\frac{N-4}{N}f_2(T_0)\quad\text{(similarly as (\ref{eq4-10}))}.\end{align}
By (\ref{eq4-10}) we obtain
\be\label{eq4-12}f_1(T_1)<\left(\frac{N-4}{N}\right)^2f_1(t_0).\ee
From (\ref{eq4-10}) and (\ref{eq4-12}), we can continue this progress to obtain a strictly increasing sequence $T_n (n\ge 0)$ such that conclusions $(1)$ and $(2)$ hold and
$$f_1(T_{2n+1})\le \left(\frac{N-4}{N}\right)^{2n+2}f_1(t_0)\quad\text{and}\quad f_2(T_{2n})\le \left(\frac{N-4}{N}\right)^{2n+1}f_1(t_0).$$
Since $T_{2n+1}$ are local minimum points of $w_1$,
Lemma \ref{lemma3-4} yields $w_1(T_{2n+1})<\la_1^\ast$. Since the function $\vp(s)=\frac{1}{2}\dd^2s^2-\frac{1}{2p}\mu_1 s^{2p}$ is strictly increasing for $s\in [0, \la_1^\ast]$,
we deduce from
$f_1(T_{2n+1})\to 0$ that $w_1(T_{2n+1})\to 0$ as $n\to\iy$. At the same time, $w_1''(T_{2n+1})\ge 0$ gives
$\bb w_2^p(T_{2n+1})<\dd^2 w_1^{2-p}(T_{2n+1})$. These imply that $w_2(T_{2n+1})\to 0$ and $T_{2n+1}\to +\iy$.
Similarly, we can prove that $w_2(T_{2n})\to 0$ and then $w_1(T_{2n})\to 0$. Combining these with conclusions $(1)-(2)$ and $T_n\to +\iy$, we easily conclude that
\be\label{eq4-13}\lim_{t\to+\iy}\min\{w_1(t), w_2(t)\}=0.\ee
This completes the proof.
\ep

We are in a position to prove the following important result, which improves Lemma \ref{lemma4-1}.

\bl\label{lemma4-2} $\lim\limits_{t\to+\iy}w_i(t)=0$ for $i=1, 2$.
\el

\bp Assume by contradiction that
\be\label{eq4-4}\limsup\limits_{t\to+\iy}w_1(t)=C_1>0.\ee

Take $s_n\uparrow +\iy$ such that $s_n$ are local maximum points of $w_1$ and $w_1(s_n)\to C_1$. By Lemma \ref{lemma4-1} we know $w_2(s_n)\to 0$. Then up to a subsequence, $w_1(\cdot+s_n)\to v_1\ge 0$ and $w_2(\cdot+s_n)\to 0$ uniformly in $C_{loc}^2(\R)$, where $v_1(0)=C_1=\max_{t\in\R} v_1(t)$ and $v_1$ satisfies $v_1''-\dd^2v_1+\mu_1v_1^{2p-1}=0$. By $K=0$ it is easy to see that
$$\frac{\dd^2}{2}v_1^2(0)-\frac{1}{2p}\mu_1v_1^{2p}(0)=0,$$
namely $C_1=\la_1$, where $\la_1$ is seen in Lemma \ref{lemma3-2}. Hence $f_1(s_n)\to 0$ as $n\to+\iy$. On the other hand, there exists $t_n<s_n$ such that $w_1'(t)>0$ for all $t\in (t_n, s_n)$ and $t_n$ is a local minimum point of $w_1$. Clearly $t_n>s_{n-1}$ and $w_1(t_n)<\la_1^\ast$ by Lemma \ref{lemma3-4}. Since
$$0<\frac{\dd^2}{2}w_1^2(t_n)-\frac{1}{2p}\mu_1w_1^{2p}(t_n)=f_1(t_n)<f_1(s_n)\to 0,$$
we obtain that $w_1(t_n)\to 0$ as $n\to+\iy$. Again, $w_1''(t_n)\ge 0$ gives
\be\label{eq4-21}\bb w_1^{p-2}(t_n)w_2^p(t_n)<\dd^2.\ee
This implies from $1<p<2$ that $w_2(t_n)\to 0$.

If there exists $T>0$ such that $w_2'(t)<0$ for all $t>T$, we may assume that $w_2'(t_n)<0$ for all $t_n$. Otherwise, $\liminf_{t\to+\iy}w_2(t)=0$ implies that $w_2$ has a sequence of local minimum points which tends to $+\iy$. Then by conclusions $(1)-(2)$ of Lemma \ref{lemma4-1}, we may also assume that $w_2'(t_n)<0$ for all $t_n$.

Similarly as the proof of Corollary \ref{corollary4-1}, up to a subsequence, $w_1(\cdot+t_n)\to 0$ and $w_2(\cdot+t_n)\to 0$ uniformly in $C_{loc}^2(\R)$. Since $w_1(s_n)\to C_1>0$, we see that $s_n-t_n\to +\iy$. By Harnack inequality, we see that
$$\frac{w_1(\cdot+t_n)}{w_1(t_n)+w_2(t_n)}\to U_1\ge 0,\quad\frac{w_2(\cdot+t_n)}{w_1(t_n)+w_2(t_n)}\to U_2\ge 0$$
uniformly in $C_{loc}^2(\R)$, where $U_i$ satisfies $U_i''-\dd^2 U_i=0$. Clearly $U_i(t)=a_i e^{\dd t}+b_i e^{-\dd t}$ for $i=1, 2$, where $a_1, a_2, b_1, b_2$ are all nonnegative constants. Since $U_1(0)+U_2(0)=1$, $U_1'(0)=0$ and $U_2'(0)\le 0$, we obtain that $a_1=b_1$, $a_2\le b_2$ and $2a_1+a_2+b_2=1$. By (\ref{eq3-2}) and $\Psi(t_n)\equiv K=0$, it is easy to get that
$$|U_1'(0)|^2+|U_2'(0)|^2=\dd^2U_1^2(0)+\dd^2U_2^2(0).$$
A direct computation gives $a_1^2+a_2b_2=0$. Hence $b_1=a_1=a_2=0$ and $b_2=1$, namely $U_1\equiv 0$ and $U_2(t)=e^{-\dd t}$. So $\frac{w_1(t_n)}{w_2(t_n)}\to 0$ and
\be\label{eq4-22}\frac{w_2(t+t_n)}{w_2(t_n)}\to U_2(t)=e^{-\dd t}\quad\text{uniformly in $C_{loc}^2(\R)$.}\ee

Now we want to apply the Harnack inequality to $\frac{w_1(t+t_n)}{w_1(t_n)}$. Remark that since $1<p<2$, it is not trivial to see that $(w_1^{p-2}w_2^{p})(\cdot+t_n)$ is locally uniformly bounded. To prove this,
we let $\eta(t)=(p-2)\ln w_1(t)+p\ln w_2(t)$, then Lemma \ref{lemma3-3} tells us that
$$|\eta'(t)|=\left|(p-2)\frac{w_1'(t)}{w_1(t)}+p\frac{w_2'(t)}{w_2(t)}\right|\le 2\dd,$$
which implies that $\eta(t)\le \eta(t_n)+2\dd |t-t_n|$. This together with (\ref{eq4-21}) yield that
\be\label{eq4-31}\bb(w_1^{p-2}w_2^{p})(t+t_n)=\bb e^{\eta(t+t_n)}\le \dd^2 e^{2\dd |t|},\quad\forall\, t\in\R.\ee
Up to a subsequence, we assume $A=\lim_{n\to+\iy}\bb w_1^{p-2}(t_n)w_2^p(t_n)$, then $A\in [0, \dd^2]$. Recall that
{\allowdisplaybreaks
\begin{align*}\frac{w_1''(\cdot+t_n)}{w_1(t_n)}-\dd^2\frac{w_1(\cdot+t_n)}{w_1(t_n)}
+\mu_1w_1^{2p-2}(\cdot+t_n)\frac{w_1(\cdot+t_n)}{w_1(t_n)}\\
+\bb (w_1^{p-2}w_2^{p})(\cdot+t_n)\frac{w_1(\cdot+t_n)}{w_1(t_n)}=0,
\end{align*}
}%
by (\ref{eq4-31}), we can apply the Harnack inequality to the above equation. Then passing to a subsequence, we deduce from (\ref{eq4-22}) that $\frac{w_1(\cdot+t_n)}{w_1(t_n)}\to V\ge 0$ uniformly in $C_{loc}^2(\R)$, where $V$ satisfies
$$V''-\dd^2 V+A e^{-p\dd t}V^{p-1}=0.$$
Clearly $V(0)=1$ and $V'(0)=0$. Recalling that $s_n-t_n\to +\iy$ and $w_1'(t)>0$ on $(t_n, s_n)$, we know that $V'(t)\ge 0$ for $t>0$ and so $V(t)\ge 1$ for all $t>0$. Define
$$H(t):=\frac{1}{2}|V'(t)|^2-\frac{\dd^2}{2}V^2(t)+\frac{1}{p}Ae^{-p\dd t}V^p(t).$$
Then $H(0)=-\frac{1}{2}\dd^2+\frac{A}{p}$.
Recalling from (\ref{eq3-7}) that $f_1(t)+ f_2(t)=\frac{\bb}{p}w_1^pw_2^p$, we have $-f_1(t+t_n)+\frac{\bb}{p}(w_1^pw_2^p)(t+t_n)>0$, namely
$$\frac{1}{2}|w_1'(t+t_n)|^2-\frac{\dd^2}{2}w_1^2(t+t_n)+
\frac{\mu_1}{2p}w_1^{2p}(t+t_n)+\frac{\bb}{p}(w_1^pw_2^p)(t+t_n)>0.$$
From this we conclude that $H(t)\ge 0$ for all $t\in\R$. In particular, $H(0)\ge 0$ gives $A\ge\frac{p}{2}\dd^2$. A direct computation gives $H'(t)=-A\dd e^{-p\dd t}V^p(t)< 0$ for $t\ge 0$, which implies that $H(t)\downarrow C\ge 0$ as $t\uparrow +\iy$. Recalling $V(t)\ge 1$ for $t\ge 0$,
we conclude that
\begin{align*}
-\frac{1}{2}\dd^2+\frac{A}{p}=H(0)=C-\int_0^{+\iy}H'(t)dt\ge A\dd\int_0^{+\iy}e^{-p\dd t}dt=\frac{A}{p},
\end{align*}
a contradiction. Therefore, (\ref{eq4-4}) is not true, namely $\lim\limits_{t\to+\iy}w_1(t)=0$. Similarly, we can prove that $\lim\limits_{t\to+\iy}w_2(t)=0$.
\ep

Now we can finish the proof of Theorem \ref{th1}.

\begin{proof}[Proof of Theorem \ref{th1} for $N\ge 5$] The fact $K(u, v)=0$ is seen in Corollary \ref{corollary4-1}. It suffices to prove that $(u, v)$ is an entire solution of (\ref{eq1}). Let $w=w_1+w_2$. It is easy to see from (\ref{eq3-1}) that $w''\ge \dd^2 w- Cw^{2^\ast-1}$. By Lemma \ref{lemma4-2}, $w''(t)>0$ for $t>0$ large and so $w'(t)\uparrow 0$ as $t\to +\iy$. That is, $w'(t)<0$ for $t>0$ large and $w(t)\downarrow 0$ as $t\to+\iy$. Then by repeating the argument of proving Theorem \ref{th1} for $N\le 4$, we obtain that $u(r)+v(r)\le C$ uniformly for $r>0$ small and so $u, v\in C^2(\RN)$. This completes the proof.
\end{proof}

Before ending this section, we would like to give an simple application of Theorem \ref{th1}. Consider the ODE system with initial values
\be\label{eq440}
\begin{cases}
\begin{split}w_{1}''-\dd^2w_{1}+\mu_1|w_1|^{2p-2}w_1+\bb |w_2|^{p}|w_1|^{p-2}w_1=0\\
w_{2}''-\dd^2w_{2}+\mu_2|w_2|^{2p-2}w_2+\bb |w_1|^{p}|w_2|^{p-2}w_2=0\\
\end{split}  \quad t\in\R,\\
w_1(0)=a_1, w_2(0)=a_2, w_1'(0)=b_1, w_2'(0)=b_2,\end{cases}\ee
where $2p=2^\ast$. Let $\Psi(t)$ be as in (\ref{eq3-2}). Then
\be\label{eq441}\Psi(t)\equiv\Psi(0)=\frac{b_1^2+b_2^2-\dd^2a_1^2-\dd^2a_2^2}{2}+
\frac{\mu_1|a_1|^{2p}+2\bb|a_1|^p|a_2|^p+\mu_2|a_2|^{2p}}{2p}.\ee
By this and the standard ODE theory, it is easy to see that (\ref{eq440}) has a solution $(w_1, w_2)$ in $\R$.
Then we have the following result.

\bt\label{th4-1}Let $(w_1, w_2)$ solve problem (\ref{eq440}) and $\Psi(0)$ be defined by initial values $a_1$, $a_2$, $b_1$, $b_2$ in (\ref{eq441}). Assume that either $\Psi(0)=0$, $a_1b_2-a_2b_1\neq 0$ or $\Psi(0)>0$.
Then at least one of $w_1, w_2$ changes sign.\et

\bp The case $\Psi(0)>0$ follows directly from Theorem B and Theorem \ref{th1}. It suffices to consider the case $\Psi(0)=0$. If neither $w_1$ nor $w_2$ changes sign, we may assume that $w_1\ge 0$ and $w_2\ge 0$. Since $a_1b_2-a_2b_1\neq 0$, we have $w_1\not\equiv 0$ and $w_2\not\equiv 0$. The maximum principle gives that $w_1, w_2>0$ in $\R$. Define
$u(x):=|x|^{-\dd}w_1(-\ln |x|)$ and $v(x):=|x|^{-\dd}w_2(-\ln |x|)$. Clearly $(u, v)$ is a solution of (\ref{eq0}) with $K(u, v)=0$. Then Corollary \ref{corollary1-1} implies that $u/v$ is a positive constant, and so $a_1b_2-a_2b_1=0$, a contradiction.\ep

\section{Non-removable singularity}
\renewcommand{\theequation}{5.\arabic{equation}}

In this section, we study singular solutions and give the proof of Theorem \ref{th2}. Comparing to Theorem \ref{th1} where the case $N\le 4$ is slightly simpler than the case $N\ge 5$, Theorem \ref{th2} is completely different. It turns out that the case $N\le 4$ seems much more complicated than the case $N\ge 5$, and here we can only deal with some special circumstances when $N\le 4$.

\begin{proof}[Proof of Theorem \ref{th2}]

{\it Step 1.} We prove the nonexistence of semi-singular solutions for $N\ge 4$.

Let $N\ge 4$. Assume by contradiction that $(u, v)$ is a semi-singular solution of (\ref{eq0}).
Without loss of generality, we assume
\be\label{eq5-1}\lim_{r\to 0}u(r)=+\iy\quad \text{and}\quad \lim_{r\to 0}v(r)=C<+\iy.\ee
By Lemma \ref{lemma3-3} and (\ref{eq5-1}), we have $\dd^{-1}|w_2'(t)|\le w_2(t)\le C e^{-\dd t}$, which implies $f_2(t)\to 0$ as $t\to+\iy$. Clearly $K(u, v)<0$.
It follows from (\ref{eq3-7}) that $f_1(t)\to-K>0$ as $t\to+\iy$. Then there exists $C_1>0$ such that
\be\label{eq5-2}\inf_{t\ge 0}w_1(t)\ge C_1\quad\text{and so}\quad u(r)\ge C_1 r^{-\dd},\quad\forall\,r\in(0, 1].\ee
In fact, if (\ref{eq5-2}) is not true, there exists $t_n\to+\iy$ such that $0\le \dd^{-1} |w_1'(t_n)|\le w_1(t_n)\to 0$, which implies $f_1(t_n)\to 0$, a contradiction.

Clearly there exists $C_2>0$ such that $v(r)\ge C_2$ for all $r\in (0, 1]$.
Recall that $v'(r)=-r^{-N/2}(w_2'(t)+\dd w_2(t))$, where $t=-\ln r$. We have $r^{N-1}v'(r)\to 0$ as $r\to 0$. Since
$$-(r^{N-1}v'(r))'=r^{N-1}(\mu_2 v^{2p-1}+\bb u^p v^{p-1}),$$
we have for $r\in (0, 1]$ that
{\allowdisplaybreaks
\begin{align*}
-r^{N-1}v'(r)&=\int_{0}^r \rho^{N-1}(\mu_2 v^{2p-1}+\bb u^p v^{p-1})d\rho\ge \bb\int_{0}^r \rho^{N-1}u^p v^{p-1}d\rho\\
&\ge C\int_{0}^r \rho^{N-1} \rho^{-\dd p}d\rho=Cr^{N/2},
\end{align*}
}%
namely $-v'(r)\ge C r^{1-\frac{N}{2}}$ for any $r\in (0, 1]$. This implies for $r\in (0, 1)$ that
$$v(r)-v(1)\ge\int_{r}^{1} C \rho^{1-\frac{N}{2}}d\rho=\begin{cases}-C\ln r,\quad\text{if $N=4$},\\
C(r^{2-\frac{N}{2}}-1),\quad\text{if $N\ge 5$}, \end{cases}$$
and so $v(r)\to +\iy$ as $r\to 0$, a contradiction with (\ref{eq5-1}). This completes the proof of Step 1.

{\it Step 2.} Let $(u, v)$ is a solution of (\ref{eq0}) with $K(u, v)<0$. We want to prove that (\ref{eq8}) holds under assumptions of Theorem \ref{th2}.

First we claim
\be\label{eq5-4}\liminf_{t\to+\iy}w_1(t)>0\quad\text{and}\quad\liminf_{t\to+\iy}w_2(t)>0. \ee

We prove this claim by discussing three cases separately.

{\bf Case 1.} $N\ge 5$.

This case is the simplest one.
Suppose that (\ref{eq5-4}) is not true.
Since Corollary \ref{corollary4-1} can not hold by our assumption $K<0$, we see that Lemma \ref{lemma3-5} does not hold and so we may assume, without loss of generality, that
\be\label{eq5-5}\liminf_{t\to+\iy}w_1(t)=C_1>0\quad\text{and}\quad\liminf_{t\to+\iy}w_2(t)=0.\ee
If $\limsup_{t\to+\iy}w_2(t)>0$, there exists a sequence of local minimum points $t_n\to+\iy$ of $w_2$ such that $w_2(t_n)\to 0$. Then the same argument as that of Corollary \ref{corollary4-1} yields $w_1(t_n)\to 0$, a contradiction with (\ref{eq5-5}). Hence $\lim_{t\to+\iy}w_2(t)=0$. It follows from (\ref{eq5-5}) that
$$w_2''(t)=w_2^{p-1}(\dd^2w_2^{2-p}-\mu_2w_2^p-\bb w_1^{p})<0\quad\text{for $t>0$ large.}$$
Since $|w_2'(t)|\le \dd w_2(t)\to 0$ as $t\to+\iy$, we see that $w_2'(t)>0$ for $t>0$ large, a contradiction with $\lim_{t\to+\iy}w_2(t)=0$. So (\ref{eq5-4}) holds.

{\bf Case 2.} $N=4$ and $\bb\ge 3\max\{\mu_1, \mu_2\}$.

Suppose that (\ref{eq5-4}) is not true. Without loss of generality, we may assume the existence of $t_n\to +\iy$ such that $w_1(t_n)\to 0$. Then up to a subsequence, $w_1(\cdot+t_n)\to 0$ and $w_2(\cdot+t_n)\to v_2\ge 0$ uniformly in $C_{loc}^2(\R)$, where $v_2$ satisfies
\be\label{eq5-6}\begin{cases}v_2''-\dd^2 v_2+\mu_2 v_2^3=0\quad\text{in $\R$},\\
\frac{1}{2}(v_2')^2-\frac{\dd^2}{2}v_2^2+\frac{\mu_2}{4}v_2^4\equiv K<0.\end{cases}\ee
Up to a subsequence as before, $\frac{w_1(\cdot+t_n)}{w_1(t_n)}\to v_1\ge 0$ uniformly in $C_{loc}^2(\R)$, where $v_1(0)=1$ and
\be\label{eq5-21}v_1''+(\bb v_2^2-\dd^2) v_1=0\quad\text{in $\R$}.\ee
The maximum principle yields that $v_1>0$ in $\R$. If $v_2$ is a constant, by (\ref{eq5-6})  we know that $v_2^2\equiv \dd^2/\mu_2$ and so $\bb v_2^2-\dd^2\equiv \frac{\bb-\mu_2}{\mu_2}\dd^2>0$, a contradiction with $v_1>0$ in $\R$. Hence $v_2$ is not a constant. Then by \cite{Fowler} it is well known that $v_2$ is a non-constant periodic function. Consequently, (\ref{eq5-21}) is the well-known Hill's equation (cf. \cite{MW}) in the literature. Remark that, although it looks very simple, Hill's equation is actually very complicated, and solutions of Hill's equation need not necessary to be periodic or bounded; see \cite{MW} for details.
Here we use the fact that $v_2'\not\equiv 0$ changes sign infinitely many times. Note that $(v_2')''+(3\mu_2 v_2^2-\dd^2) v_2'=0$. By the Sturm comparison theorem and $\bb\ge 3\mu_2$, it follows that $v_1$ must changes sign infinitely many times, also a contradiction. So (\ref{eq5-4}) holds.

{\bf Case 3.} $N=3$ and assume that (\ref{eq0}) has no semi-singular solutions.

Under our assumption, clearly
\be\label{eq6-1}\lim_{r\to 0}u(r)=\lim_{r\to 0}v(r)=+\iy.\ee
Assume by contradiction that $\liminf_{t\to+\iy}w_1(t)=0$. Recalling that $p=3$ since $N=3$, namely $w_1''-\dd^2w_1+\mu_1 w_1^5+\bb w_1^2w_2^3=0$, there exists $c>0$ such that
\be\label{eq6-2}w_1''\ge\dd^2 w_1-cw_1^2.\ee
If $\limsup_{t\to+\iy}w_1(t)>0$,
then there exists a sequence of local minimum points $t_n$ of $w_1$ such that $t_n\to+\iy$ and $w_1(t_n)\to 0$.
By (\ref{eq6-2}), there exists small $\e>0$ such that $w_1''(t)>0$ whenever $w_1(t)<2\e$.
Hence, there exist $t_n^\ast<t_n$ such that $w_1(t_n^\ast)=\e$ and $w_1'(t)<0$ for $t\in [t_n^\ast, t_n)$.
By (\ref{eq6-2}) and making $\e$ smaller if necessary, there exists a constant $c_1>0$ independent of $n$ such that
\be\label{eq6-4}t-t_n^\ast\le \frac{1}{\dd}\log\frac{w_1(t_n^\ast)}{w_1(t)}+c_1,\quad\forall\,t\in(t_n^\ast, t_n].\ee
The proof of (\ref{eq6-4}) is the same as that of (2.5) in \cite{CL3} by Chen and the second author, so we omit the details here.
Up to a subsequence, $w_1(\cdot+t_n^\ast)\to \bar{w}_1\ge 0$ and $w_2(\cdot+t_n^\ast)\to \bar{w}_2\ge 0$ uniformly in $C_{loc}^2(\R)$, where $\bar{w}_1(0)=\e$ and
$$
\begin{cases}
\bar{w}_{1}''-\dd^2\bar{w}_{1}+\mu_1\bar{w}_1^{5}+\bb \bar{w}_1^{2}\bar{w}_2^{3}=0,\\
\bar{w}_{2}''-\dd^2\bar{w}_{2}+\mu_2\bar{w}_2^{5}+\bb \bar{w}_2^2\bar{w}_1^{3}=0.\end{cases}$$
The maximum principle gives that $\bar{w}_1>0$ in $\R$.
Furthermore,
\be\label{eq6-6}\frac{1}{2}(|\bar{w}_1'|^2+|\bar{w}_2'|^2-\dd^2\bar{w}_1^2-\dd^2\bar{w}_2^2)
+\frac{1}{6}(\mu_1\bar{w}_1^{6}+2\bb \bar{w}_1^3\bar{w}_2^3+\mu_2\bar{w}_2^{6})\equiv K<0.\ee
By $|w_1'(t)|<\dd w_1(t)$ and $\dd=1/2$, we obtain that $t_n-t_n^\ast\ge 2\ln\frac{\e}{w_1(t_n)}\to +\iy$. Combining this with (\ref{eq6-4}), we easily conclude that
\be\label{eq6-5}\bar{w}_1(t)\le \e e^{\dd c_1}e^{-\dd t}=:c_2e^{-\dd t},\quad\forall\,t>0.\ee
This together with (\ref{eq6-6}) imply that $\bar{w}_2\not\equiv 0$ and then the maximum principle gives that $\bar{w}_2>0$ in $\R$. Define
$\bar{u}(x):=|x|^{-\dd}\bar{w}_1(-\ln |x|)$ and $\bar{v}(x):=|x|^{-\dd}\bar{w}_2(-\ln |x|)$. Then $(\bar{u}, \bar{v})$ is a solution of (\ref{eq0}) with $K(\bar{u}, \bar{v})=K(u, v)<0$. Moreover, (\ref{eq6-5}) implies that
$\bar{u}(x)\le c_2$ for all $|x|\in (0, 1)$, namely $(\bar{u}, \bar{v})$ is a semi-singular solution, a contradiction with our assumption. Therefore, $\lim_{t\to+\iy}w_1(t)=0$. Then there exists $T>0$ large such that $w_1(T)\le\e$ and $w_1'(t)<0$ for all $t>T$. Similarly as (\ref{eq6-4}), we see that $t-T\le \frac{1}{\dd}\log\frac{w_1(T)}{w_1(t)}+c_1$ for all $t>T$, which implies that $w_1(t)\le c_3 e^{-\dd t}$ for all $t>T$, where $c_3$ is a positive constant. In particular, $u(r)\le c_3$ for $r>0$ small, a contradiction with (\ref{eq6-1}). This finishes the proof of Case 3.

In a word, we have proved that (\ref{eq5-4}) holds. Similarly, the same conclusion as (\ref{eq5-4}) also holds for $t\to-\iy$. Hence there exists a constant $\mathcal{C}_1>0$ such that
$$\inf_{t\in\R}w_1(t)\ge \mathcal{C}_1,\quad \inf_{t\in\R}w_2(t)\ge \mathcal{C}_1.$$
Recalling that there exists $\mathcal{C}_2>0$ such that $w_1, w_2\le \mathcal{C}_2$ for all $t\in \R$, we conclude that
(\ref{eq8}) holds. This completes the proof.
\end{proof}

\end{document}